\documentstyle[epsf,sectioneqno]{article}
\newtheorem{theorem}{Theorem}[section]

\newtheorem{con}[theorem]{Conjecture}

  \title{  Absolute continuity and spectral concentration for slowly decaying potentials }
\author{    B.M. Brown, M.S.P. Eastham, D.K.R. M$^c$Cormack \\ 
 Department of Computer Science,  University of Cardiff, Cardiff, CF2 3XF, U.K. 
  }
\date{}
\begin{document}
\maketitle
\begin{abstract}
We consider the spectral function $\rho(\mu)$ $(\mu \geq 0)$
for the Sturm-Liouville equation $y^{''}+(\lambda-q)y =0$
on $[0,\infty)$ with the boundary condition $y(0)=0$ and where $q$ has slow 
decay $O(x^{-\alpha})$ $(a>0)$ as $x\rightarrow \infty$.
We develop our previous methods of locating spectral concentration for $q$ with rapid exponential decay (this Journal 81 (1997) 333-348) to deal
with the new theoretical and computational complexities
which arise for slow decay.
\newline
\vspace{0.2in}
\newline
\vspace{0.2in}
\end{abstract}
\renewcommand{\baselinestretch}{2.0}
\large \normalsize     

 \newcommand{\beq}{\begin{equation}}
\newcommand{\enq}{\end{equation}}
\newcommand{\s}{ {\cal{x}}}

\section{Introduction}
In a recent paper \cite{BEM97},   we gave a new formula
\beq
\rho^{'} (\mu)=\pi^{-1}s  \exp( -s^{-1}\int_0^\infty q(x) \sin 2 \theta(x,\mu) dx) \label{eq:1.1}
\enq
for the derivative of the spectral function associated with the Sturm-Liouville equation
\beq
y^{''}(x)+\{\lambda -q(x)\}y(x)=0\;\;\;\;(0\leq x<\infty ) \label{eq:1.2}
\enq
and the  Dirichlet boundary condition
\beq
y(0)   =0, \label{eq:1.3}
\enq
the formula applying to the situation where
\beq
q(x) \in L(0,\infty). \label{eq:1.4}
\enq
In (\ref{eq:1.1}), $\mu >0, s =\sqrt{\mu}$ and $\theta(x,\mu)$ is the solution of the first-order differential equation
\beq
\theta^{'}(x,\mu)=s-s^{-1 } q(x) \sin^2 \theta (x,\mu)
\label{eq:1.5}
\enq
such that
\beq
\theta(0,\mu)=0. \label{eq:1.6}
\enq
Also, if (\ref{eq:1.3}) is replaced by the usual general condition
\begin{displaymath}
y(0)\cos \alpha + y^{'}(0) \sin \alpha =0 \;\;\;(0 < \alpha < \pi),
\end{displaymath}
the initial factor $s$ in (\ref{eq:1.1}) is replaced by
\begin{displaymath}
s \csc^2 \alpha (s^2 +\cot^2 \alpha )^{-1}
\end{displaymath}
and (\ref{eq:1.6}) is replaced by
\begin{displaymath}
\theta(0,\mu)=-\tan^{-1}(s \tan \alpha)\;\;\;(-\pi < \theta(0,\mu)<0).
\end{displaymath}
In this paper however, we keep to (\ref{eq:1.3}) purely for simplicity.
\par
Following \cite{BEM97} and \cite{MM90} (see also \cite{BEM97a}), we say that the problem (\ref{eq:1.2})-
(\ref{eq:1.3}) exhibits {\it spectral concentration} at a point $\mu_0 \; (>0)$ if $\rho^{'}$ has a local maximum at $\mu_0$. Then $\rho$ itself has a relatively sharp increase at $\mu_0$. In \cite{BEM97}
we used (\ref{eq:1.1}) to develop computational procedures, based on (\ref{eq:1.5}) and (\ref{eq:1.6}), for locating spectral concentration points $\mu_0$ with the emphasis on potentials $q$ which decay rapidly as $x \rightarrow \infty$.
In particular, we identified a transitional property of $\theta(x,\mu)$ as $\mu$ increases through $\mu_0$ which provides a sensitive test of even slight spectral concentration \cite[sections 2-3]{BEM97}.  This property also features in  this paper, and we give details of it later at the end of section 4.
\par
A typical example considered in \cite[section 3.1]{BEM97} and suggested by \cite[example 166]{PXF91} is
\beq
q(x)=-c e^{-x/4} \cos x \;\;\;(c >0) \label{eq:1.7}
\enq
with exponential decay. However, in \cite[section 5]{BEM97}, we pointed out that our procedures are less reliable in cases of slower decay such as
\beq
q(x)=-c(1+x)^{-a} \cos x \label{eq:1.8}
\enq
with $a>1$, and we also raised the further question of what can be said when $0 < a \leq 1$, in which case (\ref{eq:1.4}) fails and it is no longer clear that (\ref{eq:1.1}) is available.
\par
In this paper, we deal with these two outstanding matters. First, in section 2, we show that (\ref{eq:1.1}) continues to hold in certain situations where (\ref{eq:1.4}) fails and $q$ is only conditionally integrable on $(0,\infty)$. 
Certain values of $\mu$ have to be avoided because of the possibility of embedded eigenvalues and discontinuities in $\rho$ but, in $\mu$ -- intervals which avoid these values, the spectrum is absolutely continuous.  Then, in section 3, we show how to extend the computational procedures in \cite{BEM97} to locate reliably spectral concentration points for  slow-decay  examples such as (\ref{eq:1.8}), with $ a \leq 1$ allowed.
\par
Finally in this introduction, we note two other approaches to spectral concentration. The software package SLEDGE
\cite{EFP96},\cite{PF93}, \cite{PXF91} replaces $q$ in (\ref{eq:1.2}) by an approximating step-function  $\tilde{q}$ over a large interval $(0,b)$. It then
computes the spectral function $\rho(\tilde{q},b,\mu)$ for (\ref{eq:1.2})
with $\tilde{q}$ and the boundary conditions (\ref{eq:1.3}) and $y(b)=0$.
When refinements of $\tilde{q}$ and $b$ lead to a stable output, the step-function $\rho(\tilde{q},b,\mu)$ provides an approximation to $\rho(\mu)$ to within a prescribed tolerance. SLEDGE is not restricted to (\ref{eq:1.4})
but, when  more explicit formulae  such as (\ref{eq:1.1}) and (\ref{eq:1.5}) are available,
these formulae provide a more sensitive means of detecting spectral concentration \cite{BEM97}.
\par
The other approach  to spectral concentration lies in the wider theoretical context of quantum resonances and spectral stability \cite{HS}.
The resonances are non-real singular spectral points and, in a future paper, we intend to develop the connection between these non-real points and the real spectral concentration points found in this paper and in \cite{BEM97}.
\section{Absolute continuity of the spectrum}
The theory from which (\ref{eq:1.1}) is derived in \cite{BEM97} was developed in the original work of Titchmarsh \cite[section 5.7]{ECT62} and Weyl \cite[p. 264]{Weyl10}.
This theory uses the fact that, subject to (\ref{eq:1.4}), there are solutions of (\ref{eq:1.1}) which  together with their derivatives are asymptotic to ${\rm \exp}(\pm i x \sqrt{\lambda})$ and $\pm i \sqrt {\lambda} {\rm \exp}(\pm  ix \sqrt{\lambda})$ as $x \rightarrow \infty$.
More recently, the existence of solutions with these  or  similar asymptotic forms, irrespective of (\ref{eq:1.4}), provides an application of the subordination theory of Gilbert, Pearson and Stolz \cite{GP87},\cite{GS92} and leads to the absolute continuity of $\rho(\mu)$ in appropriate intervals.
Such asymptotic forms are obtained by transforming (\ref{eq:1.2}) into a first-order differential system to which the Levinson asymptotic theorem \cite[section 1.3]{MSPE89} is applicable, and the necessary transformation methods were developed by Harris, Lutz and Eastham \cite{HL74}, \cite{HL75},\cite{HL77},\cite{MSPE89}. 
In particular, Behncke \cite{HB91}, \cite{HB91a}, \cite{HB94} used the transformations in \cite{HL75} to establish the absolute continuity for potentials such as (\ref{eq:1.8}) when $a >1/3$ provided that certain resonance values of $\mu$ are avoided.
In this section, we use the transformation due to Eastham and McLeod \cite[sections 4.6-4.7]{MSPE89}, \cite{EML89} to both extend this result and establish (\ref{eq:1.1})
for all $a>0$.
\par
In \cite[(4.1.8)]{MSPE89}   (\ref{eq:1.2}) is considered with $\lambda=1$ and therefore some minor changes are required in the transformation theory as presented in \cite[section 4.6]{MSPE89}. As in \cite[section 4.1]{MSPE89}, we take $q$ to have the form
\beq
q(x)=\xi(x) p(x) \label{eq:2.1}
\enq
where $p$ has period $2 \pi$ and
\begin{displaymath}
\xi(x) \rightarrow 0\;\; (x \rightarrow \infty), \;\; \xi^{'}(x) \in L(0,\infty).
\end{displaymath}
Also,
\beq
\xi(x) \not \in L^M(0,\infty), \;\; \xi(x) \in L^{M+1}(0,\infty) \label{eq:2.2}
\enq
for some integer $M \;(\geq1)$. Thus $\xi(x) = (1+x)^{-a}\;\;(0<a\leq 1)$
is the simplest example. The appropriate formulation of (\ref{eq:1.2}) as a system is
\beq
W^{'}(x)=\{i \Lambda_0 +R(x) \} W(x) \label{eq:2.3}
\enq
where, as in \cite[section 4.1]{MSPE89},
\beq
\Lambda_0={\rm dg }(\sqrt{\lambda},-\sqrt{\lambda}), \;\;R(x)=-i \xi(x) p(x)D\Omega
\label{eq:2.4}
\enq
with $D={\rm dg}(\frac{1}{2},-\frac{1}{2})$ and $\Omega$ having all entries unity.
The connection between $y$ and $W$ is
\beq
\left (
\begin{array}{c} y \\ y^{'} \end{array}
\right )
=
\left (
\begin{array}{cc} 1 &1  \\ i\sqrt{\lambda} & -i \sqrt{\lambda} \end{array}
\right ) W. 
\label{eq:2.5}
\enq
The resonance set $\sigma$ in \cite[(4.1.18)]{MSPE89} is now replaced by 
\beq
\sigma = \{ N^2/4 ; N=1,2,...\}. \label{eq:2.6}
\enq
\par
We can now proceed as in \cite[Lemma 4.6.1]{MSPE89} with a transformation
\beq
W=\{ {\rm \exp} ( \xi P_1 + \xi^2 P_2 +...+ \xi^M P_M) \} Z, \label{eq:2.7}
\enq
where the matrices $P_m$ have period $2 \pi$, ${\rm dg} P_m=0$, and
$M$ is as in (\ref{eq:2.2}).
This takes (\ref{eq:2.3}) into
\beq
Z^{'}=(i \Lambda_0 + \xi \Lambda_1 +...+\xi^M \Lambda_M + S ) Z, \label{eq:2.8}
\enq
where $S \in L(0,\infty)$ and $\Lambda_m$ is diagonal with period $2 \pi$.
Here $\lambda$ is excluded from the set $\sigma$ in (\ref{eq:2.6}).
Further, since $R$ has trace zero in (\ref{eq:2.4}), it follows that 
\begin{displaymath}
{\rm tr} \Lambda_m ={\rm -tr} P^{'}_m=0 \;\;(m\geq 1)
\end{displaymath}
as in \cite[Lemma 4.9.2 (i)]{MSPE89}.
Thus $\Lambda_m$ has the form
\begin{displaymath}
\Lambda_m = {\rm dg} (\lambda^{(m)}, -\lambda^{(m)}).
\end{displaymath}
Finally, when $\lambda$ is real and positive, the conditions of \cite[Theorem 4.6.1] {MSPE89} are satisfied by (\ref{eq:2.3}) and (\ref{eq:2.4}), and then the $\Lambda_m$ are all pure imaginary.
\par
The method used in \cite[ section 5.7]{ECT62} for obtaining the Titchmarsh-Weyl $m(\lambda)$ function and the spectral function $\rho(\mu)$ can be adapted to the situation which we have now in (\ref{eq:2.1})-(\ref{eq:2.8}).
We note that (\ref{eq:2.7}) has the form
\begin{equation}
W=(I+Q)Z, \label{eq:2.9}
\end{equation}
where $Q(x)=o(1)\;\;(x \rightarrow \infty)$, and we choose $X$ so that 
$(I+Q)^{-1}$ exists in $[X,\infty)$. We also note that (\ref{eq:2.8})
is
\begin{equation}
Z^{'} =(\Lambda + S)Z, \label{eq:2.10}
\end{equation}
where $\Lambda$  has the form
\begin{equation}
\Lambda = {\rm dg }(\nu,-\nu) \label{eq:2.11}
\end{equation}
with 
\begin{equation}
\nu =i \sqrt {\lambda} + o(1)\;\;\;(x \rightarrow \infty). \label{eq:2.12}
\end{equation}
Further, $\nu$ is pure imaginary when $\lambda$ is real  and positive.
\par
The usual integral form of (\ref{eq:2.10}) is
\begin{displaymath}
Z(x) =\Phi(x)Z(X) + \int_X^x \Phi(x)\Phi^{-1}(t)S(t)Z(t)dt
\end{displaymath}
with 
\begin{equation}
\Phi(x) = {\rm \exp } \left \{ {\rm dg } \left ( \int_X^x \nu(t)dt, - \int_X^x \nu(t) dt\right  ) \right \}. \label{eq:2.13}
\end{equation}
Then, by (\ref{eq:2.9}), the solutions of (\ref{eq:2.3}) satisfy
\begin{eqnarray}
W(x) &=& \{ I + Q(x) \} \Phi(x) \{I+Q(X)\}^{-1} W(X) \nonumber \\
& + & \{ I + Q(x) \} \int_X^x \Phi(x) \Phi^{-1}(t) \tilde{S}(t) W(t) dt,
\label{eq:2.14}
\end{eqnarray}
where
\begin{equation}
\tilde{S}=S(I+Q)^{-1} \in L(X,\infty). \label{eq:2.15}
\end{equation}
\par
If $\lambda$ and $\sqrt{\lambda}$ satisfy $0 \leq {\rm arg } \lambda < \pi$ and
$0 \leq {\rm arg } \sqrt{\lambda} < \frac{1}{2}\pi$, it follows immediately from (\ref{eq:2.12})-(\ref{eq:2.15}) and a Gronwall inequality that $W(x) {\rm \exp}
\left ( \int_X^x \nu(t) dt \right )$ is bounded on $[X, \infty )$.
To use this property in the integral term in (\ref{eq:2.14}), we define $I_1 = {\rm dg}(1,0)$ and $I_2= {\rm dg}(0,1)$ in order to split the two entries in $\Phi$ in (\ref{eq:2.13}).
Then, by (\ref{eq:2.14}), we have
\begin{displaymath}
W(x)=\{I_2 B(\lambda)+o(1) \}
{\rm \exp }\left ( - \int_X^x \nu(t) dt \right ) \;\;\; (x \rightarrow \infty)
\end{displaymath}
when $\lambda$ is non-real, where 
\begin{displaymath}
B(\lambda) = \{I+Q(X)\}^{-1} W(X)   +
\int_X^\infty {\rm \exp } \left ( \int_X^t \nu(u)du\right ) \tilde{S}(t) W(t) dt
\end{displaymath}
(cf. \cite[ (5.7.5)-(5.7.8)]{ECT62}).
Also, when  $\lambda$ has a real and positive value $\mu$, (\ref{eq:2.14}) again gives
\begin{eqnarray*}
W(x)&=&I_1 A(\mu) {\rm \exp } \left ( \int_X^x \nu(t) dt \right ) \\
&+& I_2 B(\mu){\rm \exp } \left ( -\int_X^x \nu(t) dt \right )+o(1),
\end{eqnarray*}
where $\nu$ is now pure imaginary and $A(\mu)$ is the same as $B(\mu)$
but with $-\nu$ instead of $\nu$ (cf. \cite[  (5.7.2)-(5.7.3)]{ECT62}).
Finally, in terms of the first component $a(\mu)$ of $A(\mu)$ and  the second component $b(\lambda)$ of $B(\lambda)$,
the transformation (\ref{eq:2.5}) back to $y$ gives
\begin{equation}
y(x)= \{ b(\lambda)+o(1) \} {\rm \exp } \left ( - \int_X^x \nu(t) dt \right )
\label{eq:2.16}
\end{equation}
for $\lambda$ non-real and, when $\lambda=\mu$,
\begin{eqnarray}
y(x)=a(\mu) {\rm \exp } \left ( \int_X^x \nu(t) dt \right )
+ b(\mu) {\rm \exp } \left ( -\int_X^x \nu(t) dt \right ) +o(1) \nonumber \\
y^{'}(x)=i \sqrt{\mu} a(\mu) {\rm \exp } \left ( \int_X^x \nu(t) dt \right )
-i \sqrt{\mu} b(\mu) {\rm \exp } \left ( -\int_X^x \nu(t) dt \right ) +o(1).  \nonumber \\
\label{eq:2.17}
\end{eqnarray}
We now have the same type of asymptotic formulae as in \cite[  section 5.7]{ECT62}, from which (\ref{eq:1.1}) follows as in \cite{BEM97}.
We indicate the details briefly, the only proviso being that $\mu \not \in \sigma$ in (\ref{eq:2.6}) as already mentioned.
\par
Let $y_1(x,\lambda)$ and $y_2(x,\lambda)$ be the solutions of (\ref{eq:1.2}) which satisfy the initial conditions
\begin{displaymath}
y_1(0,\lambda)=1,\;\;\; y_1^{'}(0,\lambda)=0,\;\;\; y_2(0,\lambda)=0, \;\;\;y_2^{'}(0,\lambda)=1,
\end{displaymath}
and let $a_1,b_1,a_2,b_2$ denote the corresponding multipliers as in (\ref{eq:2.16}) and (\ref{eq:2.17}).
Then \mbox {$y_1 + m y_2 \in L^2(0,\infty)$} gives
\begin{equation}
m(\lambda)=- b_1(\lambda)/b_2(\lambda) \;\;\;\;({\rm im } \lambda \neq 0)
\label{eq:2.18}
\end{equation}
as in \cite[(5.7.9)]{ECT62}.
When  $\lambda = \mu$, $y_1$ and $y_2$ are real-valued and hence
\begin{displaymath}
a_j(\mu)= \overline { b_j(\mu)}\;\;(j=1,2)
\end{displaymath}
in (\ref{eq:2.17}).  Then $W(y_1,y_2)=1$ gives
\begin{displaymath}
{\rm im }(\overline{b_2} b_1)(\mu) = - ( 4 \sqrt{\mu})^{-1}.
\end{displaymath}
Hence, by (\ref{eq:2.18}),
\begin{eqnarray*}
\lim_{\lambda \rightarrow \mu } {\rm im } \;m(\lambda) &=&  
\{ 4 \sqrt{\mu} \mid b_2 (\mu) \mid ^2 \}^{-1} \\
&=&\lim_{x \rightarrow \infty }
\{ \mu^{1/2} y_2^2(x,\mu)+\mu^{-1/2}y_2^{'2}(x,\mu) \}^{-1}
\end{eqnarray*}
by (\ref{eq:2.17}) and (\ref{eq:1.1}) follows as in \cite{BEM97}.
\par
Finally in this section, we note that the requirement $\mu \not \in \sigma$ can be relaxed as follows when $M=1$ in (\ref{eq:2.2}).
Let $c_n \;(-\infty < n <\infty)$ denote the complex Fourier coefficients of $p(x)$ in (\ref{eq:2.1}). Then it is shown in \cite[section 4.2]{MSPE89} ( see also \cite{MSPE91})
that (\ref{eq:2.17}) continuous to hold when $\mu = \frac{1}{4}N^2$ for some $N$ provided that $c_N=0$. Thus, altogether, (\ref{eq:1.1}) is valid
\begin{enumerate}
\item
for all $\mu >0$ when $\xi \in L(0,\infty)$,
\item
for all $\mu >0$ except those $\mu=\frac{1}{4}N^2$ for which $c_N \neq 0$, when $M=1$ in (\ref{eq:2.2}),
  \item
for  
$ \mu >0 $ and $\mu  \not \in \sigma$ when $ M \geq 2$   in(\ref{eq:2.2}).
 \end{enumerate}
\section{An integration algorithm}
We aim to compute $\rho^{'}$ in (\ref{eq:1.1}) to within a reasonable degree of accuracy such as $10^{-5}$. The error in $\rho^{'}$ is of course made up of a number of components: the truncation error due to the approximation of the semi infinite interval by a finite interval,   the error inherent in the  solving 
algorithm of the differential equation  and the  rounding error due to the rational approximation of real numbers in computer arithmetic.
We focus here on the truncation error and rely on the standard theory for both the error in the Runge Kutta algorithm as well as the floating-point numerical software on our computer system.
The infinite integral in (\ref{eq:1.1}) is of course truncated at a suitable value $X_0$.
However, to achieve a truncation error of $10^{-6}$ in the integrand when, for example,
$a=3$ in (\ref{eq:1.8}) requires $X_0=10^3$, and integration over the large range $(0,10^3) $ is unreliable.
The situation is much worse for smaller values of $a$.
In this section, we develop an iterative algorithm which accelerates the convergence of the integral in (\ref{eq:1.1}) and enables us to cope with potentials such as (\ref{eq:1.8}) when $a >0$.
\par
Guided by the example (\ref{eq:1.8}), we give the algorithm for the case $p(x)=\cos x$, so that (\ref{eq:2.1}) is
\begin{equation}
q(x) = \xi (x) \cos x \label{eq:3.1}
\end{equation}
and $\xi(x)$ is as before but with $M\geq 0$ in (\ref{eq:2.2}).
Our methods also cover the more general situation where $p(x)$ is a finite
Fourier series, but the details become more complicated.
We require the following trigonometric identity, valid for any $\alpha, \beta, \theta$ and $x$:
\begin{eqnarray}
 && 8 \cos x \sin^2 \theta \sin ( \alpha \theta + \beta x ) \nonumber \\
&=& 2 \sin \{ \alpha \theta + (\beta +1) x \} + 2 \sin \{ \alpha \theta + (\beta -1) x \} \nonumber \\
&-& \sin \{ (\alpha +2)\theta + (\beta +1 ) x \} -\sin \{ (\alpha -2 ) \theta + 
(\beta -1) x \} \nonumber \\
&-& \sin \{ (\alpha +2)\theta + (\beta -1 ) x \} -\sin \{ (\alpha -2 ) \theta + 
(\beta +1) x \}. \label{eq:3.2}
\end{eqnarray}
This is easily verified, and there is a similar identity with $\cos (\alpha \theta + \beta x)$ on the left and all cosines on the right.
\par
Next we require the following integrals over $[0,\infty)$ with
$\theta$ as in (\ref{eq:1.1}), $ F \in L(0,\infty)$ and $F(\infty)=0$:
\begin{eqnarray}
I(F,\alpha,\beta) &=& \int F(x) \sin (\alpha \theta + \beta x ) dx \nonumber \\
J(F,\alpha,\beta) &=& \int F(x) \cos (\alpha \theta + \beta x ) dx \nonumber \\
K(F,\alpha,\beta) &=& \int F(x)\cos x \sin^2 \theta \sin (\alpha \theta + \beta x ) dx \nonumber \\
L(F,\alpha,\beta) &=& \int F(x) \cos x \sin^2 \theta\cos (\alpha \theta + \beta x ) dx. \label{eq:3.3}
\end{eqnarray}
It follows from (\ref{eq:1.5}), (\ref{eq:3.1}) and an integration by parts that 
\begin{eqnarray}
&& (\alpha s + \beta) I ( F, \alpha, \beta) \nonumber \\
&=& \int F \sin (\alpha \theta + \beta x ) ( \alpha \theta^{'} + \beta + \alpha s^{-1} \xi \cos x \sin^2 \theta ) dx \nonumber \\
&=& F(0) + J(F^{'},\alpha,\beta ) + \alpha s^{-1}K(F\xi,\alpha,\beta)
\label{eq:3.4}
\end{eqnarray}
and similarly
\begin{equation}
 (\alpha s + \beta) J ( F, \alpha, \beta)=-I(F^{'},\alpha,\beta) + \alpha s^{-1} L (F\xi,\alpha,\beta). \label{eq:3.5}
\end{equation}
Finally in these introductory formulae, it follows from (\ref{eq:3.2}) and (\ref{eq:3.3}) that
\begin{eqnarray}
8K(F,\alpha,\beta) &=& 2 I(F,\alpha,\beta+1)+2 I(F,\alpha,\beta-1) \nonumber \\
&-& I(F,\alpha+2,\beta +1)-I(F,\alpha-2,\beta-1) \nonumber \\
&-&I(F,\alpha+2,\beta-1)-I(F,\alpha-2,\beta+1) \label{eq:3.6}
\end{eqnarray}
\begin{eqnarray}
8L(F,\alpha,\beta) &=& 2 J(F,\alpha,\beta+1)+2 J(F,\alpha,\beta-1) \nonumber \\
&-& J(F,\alpha+2,\beta +1)-J(F,\alpha-2,\beta-1) \nonumber \\
&-&J(F,\alpha+2,\beta-1)-J(F,\alpha-2,\beta+1). \label{eq:3.7}
\end{eqnarray}
\par
We can now return to the integral in (\ref{eq:1.1}), which we denote by $I_0$. Then, by (\ref{eq:3.1}), we can write
\begin{equation}
2 I_0 = I(\xi,2,1)+I(\xi,2,-1). \label{eq:3.8}
\end{equation}
By (\ref{eq:3.4}) and (\ref{eq:3.6}),  the two $I-$integrals here can be expressed in terms of $I$ and $J$ integrals with integrands containing $\xi^{'}$ and $\xi^2$.
These last integrals converge more rapidly than those
  in (\ref{eq:3.8}) for cases such as 
\begin{equation}
\xi(x) = ( {\rm const. }) (1+x)^{-a}. \label {eq:3.9}
\end{equation}
Repetition of the algorithm (\ref{eq:3.4})-(\ref{eq:3.7}) accelerates the convergence by introducing integrands with higher derivatives and higher powers of $\xi$.
\par
Certain values of $\mu$ have to be excluded to avoid a zero factor $\alpha s + \beta$ on the left-hand side of (\ref{eq:3.4}) and (\ref{eq:3.5}).
Thus, with $\alpha=2$ and $\beta=-1$ in (\ref{eq:3.8}), we exclude
\begin{equation}
\mu=1/4 \label{eq:3.10}
\end{equation}
at the first implementation of the algorithm.
At the next application of the formulae (\ref{eq:3.4}) and  (\ref{eq:3.5}) the value 
\beq
\mu=1 \label{eq:3.11a}
\enq
is also excluded and, at the third application of   (\ref{eq:3.4}) and  (\ref{eq:3.5}), the additional values
\beq
1/36,\;\;\;1/16,\;\;\;9/16,\;\;\;9/4. \label{eq:3.12a}
\enq
We note that there is an overlap of these values and the $ N^2/4$ values discussed at the end of section 2 for the validity of (\ref{eq:1.1}) when $M \geq 1$.
\section{Implementing the algorithm}
In this section we show how  the integration   algorithm from section 3 is used to identify points of spectral concentration for potentials of the form 
\begin{equation}
q(x)=-c  (1+x)^{- a }\cos x \label{eq:4.1}
\end{equation}
given by (\ref{eq:3.1}) and (\ref{eq:3.9}),     where $c>0$ and $ a >0$. 
The algorithm is implemented  using both symbolic methods and numerical approximations.
\par
We recall from (\ref{eq:3.8}) that the integral  which appears in the formula
for $\rho^{'}$ is  the sum of two integrals $I( \xi,2,1)$ and $I(\xi,2,-1)$. The first part of the algorithm  consists of a procedure to improve the convergence of these integrals. We focus first on the case $a\geq 2$ and  then comment  on the procedure that we have been forced to adopt for smaller values of $a$. 
\par
 First  an acceptable order of convergence $\epsilon$ is  decided upon: we have chosen $\epsilon = x^{-6}$. A purpose written Mathematica code is used to repeatedly  apply the integration by parts formulae   (\ref{eq:3.4})-(\ref{eq:3.7}) to the integrals whose integrands are larger then $\epsilon$.  Starting with the   integrals $I( \xi,2,1)$ and $I(\xi,2,-1)$, the integration by parts formulae generate integrals of the types $I,J,K,L$ (cf. (\ref{eq:3.3})) whose integrands have smaller order than the integrand in (\ref{eq:1.1}) together with terms that do not involve the variable $x$.  This procedure
is repeated until all integrands have order  less than or equal to $\epsilon$.  Thus Mathematica is used to generate a symbolic formula which consists of terms that do not depend upon $x$, denoted by $C(\mu)$,  together with  a sum of integrals of type (\ref{eq:3.3}).
Next  a Mathematica code is written to parse the formula and reconstruct the  integrands.  The symbolic formula is finally converted into a Fortran 77 function which for convenience we   denote by $F(x,\mu)$.
\par
The next task is to  evaluate the integral (\ref{eq:1.1}) numerically from the improved integrand for each value of $\mu$ under consideration.
As the integral $I_0$    depends upon $\theta$, which itself is a solution of  the differential equation (\ref{eq:1.5}), this is done by solving the system
\beq
\left ( \begin{array}{c} I_0(x,\mu) \\ \theta (x,\mu)\end{array} \right )^{'}=
\left ( \begin{array}{c} F(x,\mu) \\ s-s^{-1} q(x) \sin^2 \theta \end{array} \right ) \label{eq:4.2a} \enq
over $[0,X]$ for some large $X$, where $I_0(x,\mu)$ is the integral over $(0,x)$ in (\ref{eq:1.1}), subject to the initial conditions
\begin{displaymath}
I_0(0,\mu)=C(\mu),\;\;\; \theta(0,\mu)=0,
\end{displaymath}
typically X=100 has been used in our example.
\par
There are several  practical difficulties encountered in
 performing the above tasks. First the improvements in the convergence of the integrals brought about by the integration by parts   generate a large number of integrals $I,J,K,L$ as indicated by the right-hand sides of (\ref{eq:3.6}) and (\ref{eq:3.7}). Numerical inaccuracies  do not allow us at the moment to improve the integrands  beyond a certain point. The inaccuracies arise from the $K$ and $L$ type integrals when we attempt to improve the integrand beyond $\xi^3$.
The combined integrand $F(x,\mu)$ therefore involves 
\beq
\xi^{''''},\;\;\xi \xi^{'''},\;\; \xi \xi^{''}, \;\; \xi^2 \xi^{'}, \;\; \xi^{'2}, \;\; \xi^{3}
\label{eq:4.3a}
\enq
and the excluded values of $\mu$ are $\frac{1}{4}$ and $1$ as stated in (\ref{eq:3.10}) and (\ref{eq:3.11a}).
 We have also tried to  perform the numerical integration by extending the system of differential equations (\ref{eq:4.2a}) to one in which each entry is  only one of the integrands of $I,J,K,L$ together with the defining equation for $\theta$.
However this procedure, which has required a more sophisticated parsing routine to be written, has produced no significant improvement in the results.
\par
When $a \geq 2$, the order of convergence $\epsilon=x^{-6}$ is achieved
by the $F(x,\mu)$ indicated by (\ref{eq:4.3a}).
When $a <2$, we continue to use the same $F(x,\mu)$ but with a consequent increase in the truncation error.  However, the increased error leads to two complications which become more serious as $a$ decreases.
The first complication is that (\ref{eq:1.1}) only gives the approximate location of the spectral concentration points.
The second and more serious complication is that spurious maxima of $\rho^{'}$ are produced, the more so as $a$ decreases, and therefore it is necessary to identify the true maxima.
We resolve these difficulties by using a transitional property of $\theta(x,\mu)$ which we now describe.
\par
As reported in \cite[section 2]{BEM97a} ( see also \cite{BEM97}),
spectral concentration at a point $\mu_0$ is indicated by a certain transitional behaviour of $\theta(x,\mu)$ as $\mu$ increases through $\mu_0$.
Let $(x_1,x_2)$ be an interval in which the 
   Sturm-Liouville coefficient $\mu-q(x) <0$ and $\mu^{'}$ and $\mu^{''}$ be  suitably close to a point of $\mu_0$ of spectral concentration  with $\mu^{'} < \mu_0 < \mu^{''}$.
 Then we expect the integral in (\ref{eq:1.1}) to be large and negative--
and therefore producing spectral concentration-- if
\begin{equation}
(N+\frac{1}{2})\pi< \theta (x,\mu_0) < (N+1)\pi \label{eq:4.2}
\end{equation}
in $(x_1,x_2)$, where $N(\geq 0)$ is an integer.
We find that, particularly in situations of sharp concentration, (\ref{eq:4.2}) is realised with the following features.
\begin{enumerate}
\item
$\theta(x,\mu^{'})$ and $\theta(x,\mu^{''})$  are close together for $0 \leq x \leq x_1$ with their values at $x_1$ close to $(N+\frac{1}{2})\pi$.
\item
$\theta(x,\mu^{''})-\theta(x,\mu^{'})$ is close to $\pi$ for $x$ to the right of $x_2$, with   $\theta(x,\mu^{''})$ close to $(N+\frac{3}{2})\pi$.
\item
$\theta(x,\mu_0)$ is close to $(N+1)\pi$.
\end{enumerate}
Thus the graph of $\theta(x,\mu)$ undergoes a rapid transition in $(x_1,x_2)$ as $\mu$
increases from $\mu^{'}$ to $\mu^{''}$.
The transitional behaviour is illustrated by the graphs in the next section.
\par
We return now to $F(x,\mu)$ in (\ref{eq:4.3a}) and the case $a < 2$.
The first complication, concerning the approximate location of spectral concentration points, is resolved by using this approximate location as the starting point of a search range for $\mu$ within which (\ref{eq:1.5}) is solved numerically for
$\theta(x,\mu)$. Then the points $\mu_0$ at which the transition occurs can be located more precisely.
We can also use the transition property of $\theta(x,\mu)$ even when $a \geq 2$
to verify independently that true spectral concentration points obtained precisely as the local maxima of $\rho^{'}$ are not artifacts of our methods.
At the same time, any apparent maxima which are not associated with the transition
property are rejected as spurious.
We find that, as $a$ decreases, the number of spurious maxima increases significantly making the task of rejection a major (but unseen) part of our work.
\par
Finally in these comments, we note that there is a further complication as a result of the excluded values 
$\mu=\frac{1}{4}, 1,\frac{1}{36},...$ in (\ref{eq:3.10})-(\ref{eq:3.12a}).
Although is is clear from the integration by parts formulae (\ref{eq:3.4})-(\ref{eq:3.7}) that our methods of evaluating $\rho^{'}$
must exclude these points, the numerical realisation of our algorithm produces unreliable results in neighbourhoods of these points.
Again we test these neighbourhoods for the appearance of the transitional behaviour of $\theta(x,\mu)$.
We refer also to \cite[section 6]{BEM97a} for a similar use of $\theta(x,\mu)$
in a different but related situation where a direct formula for $\rho^{'}$
poses difficulties.
\par
In the case of (\ref{eq:4.1}), the intervals $(x_1,x_2)$ in (\ref{eq:4.2}) are approximately 
\beq
((2r+\frac{1}{2})\pi, (2r + \frac{3}{2}) \pi)\;(r=0,1,2,...).\label{eq:4.5a}
\enq
 We then denote $\mu_0$ in (\ref{eq:4.2})
by
\begin{displaymath}
\mu(c,N)\;(r=0), \;\nu(c,N) \;(r=1),\; \xi(c,N) \;(r=2)
\end{displaymath}
and, in the next section, we record our findings concerning the location of these points
when $c(>0)$ is regarded as a parameter. We shall also comment on higher values of
$r$ as appropriate.

\section{Examples}

\subsection{Example : $a=2$ in (\ref{eq:4.1})}
Here $q(x)$ is $L(0,\infty)$ and (\ref{eq:1.1}) is valid for all $\mu >0$.
We apply the algorithm   as described in (\ref{eq:4.3a}) and  we have to exclude values of $\mu$ near to $\frac{1}{4}$ and $1$ in the computation
of $\rho^{'}(\mu)$.
In Table 1 we list the spectral concentration points which we have located as giving
 local maxima of $\rho^{'}(\mu)$ except that, as mentioned in section  4,
we  have located points near to $\frac{1}{4}$ and $1$ by identifying the value of $\mu$ which is associated
 with the transitional property of $\theta(x,\mu)$.
Thus, for example, Figure 1 gives the graphs of $\theta(x,\mu)$ which identify
 $\mu(49.26,1)=0.25$ and provide part of the evolution of $\mu(c,1)$ as $c$ varies which is summarised in the second column of Table 1.
We point out that the graphs represent $\theta$ $( {\rm mod} \pi)$,
whence the repeated cut-offs at the ordinate $\pi$.
\par
When $c$ has the particular value $122.1$, we find that 
\begin{displaymath}
\mu(c,2)= \nu(c,3)=0.51.
\end{displaymath}
Thus the two local maxima of $\rho^{'}$ arising from $\mu(c,2)$ and $\nu(c,3)$ coalesce or, equivalently, two different intervals $(x_1,x_2)$ ( with $r=0$ and $r=1$) simultaneously make relatively large contributions to the integration in (\ref{eq:1.1}).
 One result of this coalescence is that $\nu(c,3)$ becomes $\nu(c,4)$ when $c > 122.1$. Thus, the last-line entry $0.49$ in Table 1 is in fact $\nu(125,4)$.

\begin{figure}[htbp]
\centerline{
\epsfxsize=12cm
\epsffile {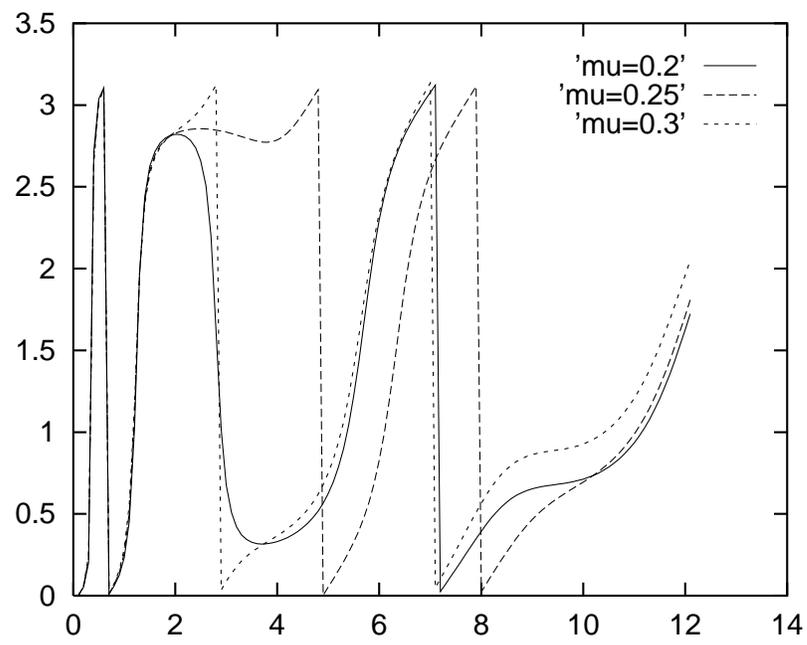 }
}
\caption{$\theta$ graph for $\mu(49.26,1)=0.2 5$ }
\label{fig:3}
\end{figure}

\begin{table}[htbp]
 
\begin{center}
\begin{tabular}{||c|c|c|c|c|c|} \hline
 $c$ & $\mu(c,0)$ & $\mu(c,1)$  & $\mu(c,2)$&$\mu(c,3)$  & $\nu(c,3)$ \\ \hline
2   &  0.45 & &  & &  \\ \hline
 7  & 0.26  & 2.05&   & & \\ \hline
9   & 0.08  &2.41 &   & & \\ \hline
20   &   & 2.10&  & &  \\ \hline
30   &   &1.87 & 5.28  & & \\ \hline
40   &   &1.29 & 5.51  & & \\ \hline
50   &   &0.15 & 5.10 & &0.68 \\ \hline
60   &   & & 4.81 & &0.72 \\ \hline
70   &   & & 4.73 & 9.75&0.70  \\ \hline
80   &   & & 4.45 & 9.95& 0.67 \\ \hline
100   &   & & 3.11 & 9.26& 0.60 \\ \hline
120   &   & & 0.80 &8.9 & 0.52 \\ \hline
125   &   & &0.10 & 8.84& 0.49   \\ \hline

   \end{tabular}
\end{center}
\caption{ $a=2$ }
\label{tab:1}
\end{table}

\newpage

\subsection{Example: $a=1$ in (4.1)}
Now $M=1$ at the end of section 2, and (\ref{eq:1.1}) is valid for all $\mu >0$ except $\mu=1/4$.
We have followed the procedure for $a<2$ as described in section 4, in which the excluded values are $\mu=1/4$ (again) and $\mu=1$.
We find that spectral concentration points occur in greater profusion than for $a=2$, and our results for $c \leq 100$ are given in Tables 2 and 3. Coalescing points occur as follows:
\begin{eqnarray*}
\mu(c,1) = \nu(c,3) = 1.98 ,\;(c=24.25) \\ 
\mu(c,1) = \xi(c,4) = 0.53,\;(c=31.2) \\
\mu(c,2) = \nu(c,5) = 4.02 ,\;(c=64.6) \\
 \end{eqnarray*}
as well as 
\begin{displaymath}
\mu(c,3) = \nu(c,8) =10.06 ,\;(c=106).  
\end{displaymath}
As  $c$ increases through these respective values, the oscillation number $N$ for $\nu(c,N)$ or $\xi(c,N)$ increases by one,
because of the transition $\pi$ in the values of $\theta$ caused by the $\mu$ spectral concentration points.
This increase is indicated at the tops of columns in Table 3.
Also, in Figure 2, we give the $\theta-$graphs which show the transitional behaviour twice  to illustrate the coalescence of $\mu(c,2)$ and $\nu(c,5)$ when $c=64.6$.
 \begin{figure}[htbp]
\centerline{
\epsfxsize=11cm
\epsffile {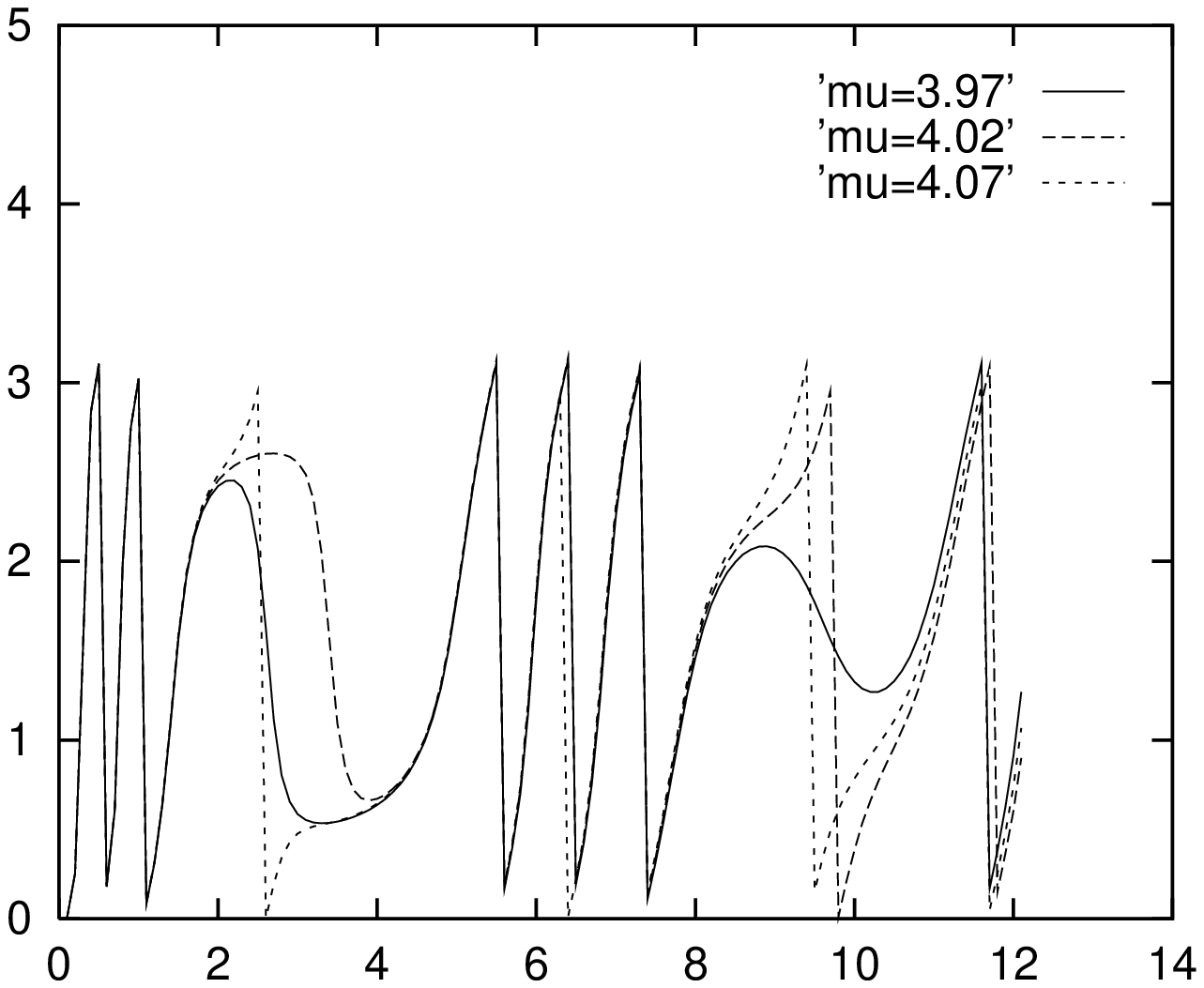}
}
\caption{$\theta$ graphs for $\mu(64.6,2)=\nu(64.6,5)=4.02$ }
\label{fig:3a}
\end{figure}

\begin{table}[hp]

\begin{center}
\begin{tabular}{||c|c|c|c|c|c|} \hline
 $c$ & $\mu(c,0)$ & $\mu(c,1)$  & $\mu(c,2)$&$\mu(c,3)$  & $\mu(c,4)$\\ \hline
2   &  0.5 & &  & & \\ \hline
 3  & 0.48 & 3.05&  & &  \\ \hline
4   & 0.41  &2.73 &  & &  \\ \hline
5   & 0.26  & 2.68&  & &  \\ \hline
6   &  0.03 &2.68 &   & &  \\ \hline
10   &   &2.76 & 6.77 & &  \\ \hline
15   &   &2.97 & 6.81 & &  \\ \hline
20   &   & 2.61& 6.84 & &   \\ \hline
25  &   &1.85 & 7.09 & 12.51&   \\ \hline
30   &   & 0.80& 7.36 & 12.73&   \\ \hline
40      & &   & 7.35&12.82&18.38 \\ \hline
50      & &   &6.43 &13.36& 20.41\\ \hline
60      & &  &4.89& 13.90&20.55  \\ \hline
70      & &  & 2.89& 13.92&20.79  \\ \hline
80      & &  & 1.75& 13.36&21.39  \\ \hline
90      & &  & 1.55& 12.35&22.06 \\ \hline
100    &  & &1.31 & 11.01 &22.47  \\ \hline

   \end{tabular}
\end{center}
\caption{ $a=1$ }
\label{tab:2}
\end{table}

\begin{table}[hp]
\begin{center}
\begin{tabular}{||c|c|c|c|c|c|} \hline
 $c$ & $\nu(c,2)$ & $\nu(c,3/4)$  & $\nu(c,5/6)$&$\nu(c,8)$  & $\xi(c,4/5)$\\ \hline
2   &    & &  & & \\ \hline
 3  &   &  &  & &  \\ \hline
4   &    &  &  & &  \\ \hline
5   & 0.70  &   &  & &  \\ \hline
6   &  0.75 &  &   & &  \\ \hline
10   & 0.72  &  &   & &  \\ \hline
15   & 0.59  &  &   & &  \\ \hline
20   & 0.41  & 2.22&   & & 0.69 \\ \hline
25  &  0.14 &2.06 &   &  & 0.64  \\ \hline
30   &   & 2.05& 4.05 &  & 0.57  \\ \hline
40      & & 1.83  & 4.10& &0.34 \\ \hline
50      & & 1.46  &4.15 & & 0.19\\ \hline
60      & &1.00  &4.08&  &0.04  \\ \hline
70      & & 0.43 & 3.88& &  \\ \hline
80      & &  & 3.60& &  \\ \hline
90      & &  & 3.24& & \\ \hline
100    &  & &2.82 & 10.14 &   \\ \hline

   \end{tabular}
\end{center}
\caption{ $a=1$ }
\label{tab32}
\end{table}
 \newpage
\subsection{Example : $a=\frac{1}{2}$ in (4.1)}
We have intimated in sections 2 and 4 that slow decay spawns theoretical and computational complexities, the latter including the task of segregating
the large number of spurious and actual local maxima of $\rho^{'}$.
This task is expensive in computer time  and therefore we have restricted
the range of $c$ in this part of our investigation to $0 < c\leq 30$.
Our findings for $\mu, \nu$ and $ \xi$ points are summarised in
Table 4, in which figures are given to more than two decimal places when necessary to distinguish between closely situated points.
\par
In particular, we find that there are two values of $c$ and $\mu$ in whose  neighbourhoods several very close (but apparently not coalescing)
points of spectral concentration exist. The values are
\begin{displaymath}
c=3,\;\; \mu=0.49; \;\; \; \; c=5.5,\;\; \mu=0.66.  
\end{displaymath}
In the case of $c=3$, for example, we have evidence for the existence of
at least seven such spectral concentration points, corresponding to $0 \leq r \leq 6$ in (\ref{eq:4.5a}). In Figure 3,  we exhibit this evidence for
$0 \leq r \leq 4 $.
The figure shows $\theta-$graphs for a decreasing sequence of 
values of $\mu$ (all    close to $0.49$) where the transition occurs for
$r=0,1,2,3,4$ in turn.
We make a further comment on this matter  in section 6.3 below.
\par
For $c$ in the stated range $(0,30)$ we have identified  one pair of coalescing
points which, as in section 5.2, are associated with a change in the oscillation
number $N$, as follows: $\mu(c,1)=\nu(c,3)=0.67\;\;\;(c=24.5)$.
Further, in addition to what is reported in Table 4, we have found that, for
$15 \leq c \leq 24.5$,  there is a $\nu(c,4)$ which is almost identical to $\mu(c,1)$.
Thus, for $c$ in this range, we have a spectral concentration point which enjoys
contributions from both $r=0$ and $r=1$ in (\ref{eq:4.5a}).
This linkage between $\mu$ and $\nu$ points is a new feature for $a=1/2$
and it poses the question whether there is a theoretical explanation.
\footnotesize
\begin{table}[hp]

\begin{center}
\begin{tabular}{||c|c|c|c|c|c|c|c|c|c|c|} \hline
 $c$ & $\mu(c,0)$ & $\mu(c,1)$  & $\mu(c,2)$&$\nu(c,1)$  &$\nu(c,2)$&$\nu(c,3\backslash 4)$  &$\xi(c,2)$&$\xi(c,3)$&$\xi(c,4)$&$\xi(c,5)$\\ \hline
2 &0.61 & & &0.56& & & & & &     \\ \hline
 3  & 0.50 & & & 0.487& 0.73&  &0.4869 &0.62 & &     \\ \hline
4   & 0.26 & & & &0.72  & & &0.66 &0.78 &       \\ \hline
 6  &  & & & & 0.575&1.83  & &0.569 &0.71 &       \\ \hline
8 &  &3.49 & & &0.33 &2.08  & &0.331  &0.58 &1.85   \\ \hline
10  &  &3.61 & & &0.09 &2.11  & & & &  2.01 \\ \hline
15  &  &3.11 & & & & 1.842   & & & &  1.841      \\ \hline
20  &  & 2.0& 8.82 & & &1.2986    & & & & 1.2984       \\ \hline
25  &  &0.51  &9.02 & & &0.59 &
 & & &        \\ \hline
30  &  &       &8.72 & & & & & & &        \\ \hline
    \end{tabular}
\end{center}
\caption{  $a=1/2$ }
\label{tab:2a}
\end{table}

\newpage
\normalsize
\begin{figure}[htbp]
\centerline{
\epsfxsize=16cm
\epsffile {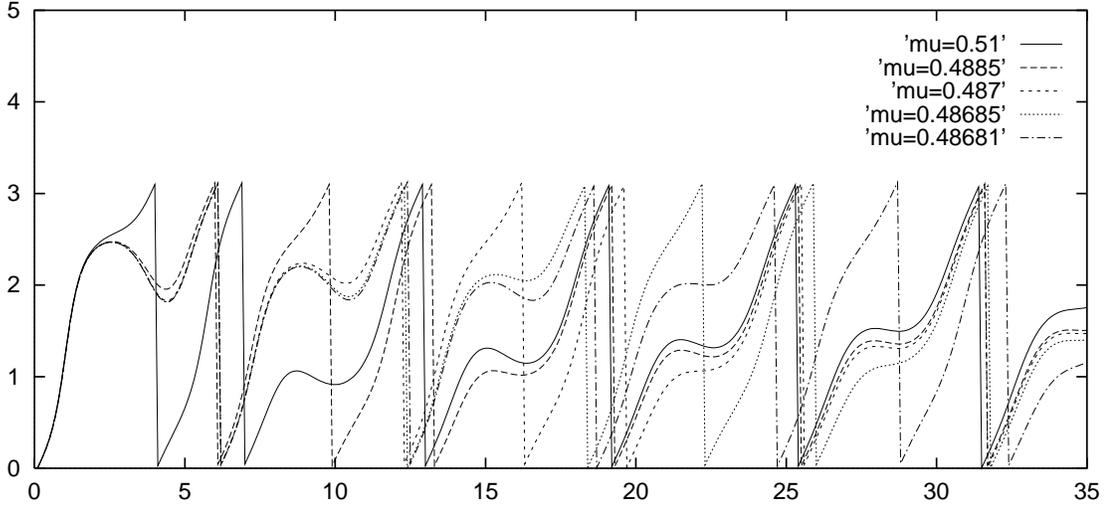}
}
\caption{$\theta$ graphs for $a = \frac{1}{2} , c = 3$ }
\label{fig:4}
\end{figure}

\section{Concluding remarks}
\subsection{Absence of spectral concentration}
In \cite[section 5]{BEM97}, we made the following conjecture concerning the appearance of spectral concentration for (\ref{eq:1.2}) and (\ref{eq:1.3}) when $q(x)$ has the form $cQ(x)\;\;(c>0)$ such as we have in (\ref{eq:1.7}),
(\ref{eq:1.8}) and (\ref{eq:4.1}).
\begin{con}
Let $Q(x)<0$ in some interval $(0,x_0)$.
Let $Q(x)$ change sign one or more times as $x$ increases with $Q(x)$ finally decaying to zero as $x\rightarrow \infty$ and $Q(x) \in L(0,\infty)$.
Then there is a number $c_0\;(>0)$ such that spectral concentration does not appear for any $\mu>0$ when $c$ lies in the range $(0,c_0)$. 
\end{con} 
\par
For the example (\ref{eq:1.7}) we showed in \cite{BEM97} that the conjecture is true  with $0.28 < c_0 \leq 0.29$ and, more generally, an affirmative answer to the conjecture has been given recently for the special case of (\ref{eq:2.1})
in which $x \xi (x)$ is $L(0,\infty)$.
Thus (\ref{eq:4.1}) is covered provided that $ a >2$.
The conjecture remains undecided when $1 < a \leq 2$ in (\ref{eq:4.1}), 
although our computational findings for $a=2$ indicate that $c_0$ exists with $1.0 < c_0 \leq 1.1$.

 \subsection{More general potentials}
There is no difficulty in principle in extending our
methods for (\ref{eq:2.1}) to potentials which are a finite sum
\beq
q(x)=\sum_1^L \xi_l (x) p_l (x) \label{eq:6.1}
\enq
and the $p_l(x)$ have different periods $\omega_l$,
provided that the $\omega_l$ are mutually commensurable.
The $p_l(x)$ would then all have a common basic period $\omega$,
and the set $\sigma$ in (\ref{eq:2.6})
is modified to
\beq
\sigma = \{ N^2 \pi^2 /\omega^2; N=1,2,...\}. \label{eq:6.2}
\enq
However, in the absence of a special form such as (\ref{eq:2.1}) or (\ref{eq:6.1}), it is no longer clear what can be said in general about the absolutely continuous nature of $\rho$ when $q$ fails to be $L(0,\infty)$.
There is no simple exceptional set $\sigma$ such as (\ref{eq:6.2}) because in 
\cite[Section 4.4]{EK}, \cite{EML77} and \cite{TE80} it is shown that, given any set of isolated  positive real numbers $\mu_n$ and given $a\;(0<a<1)$,
there is a potential $q(x)=O(x^{-a})\;(x \rightarrow \infty)$, such that the spectral function has discontinuities at the $\mu_n$.
\subsection{Higher values of $r$ in (4.5)}
We have focused on $r=0,1,2$ in (\ref{eq:4.5a}) in order to establish our computational methods, and consequently we have located  mainly $\mu,\nu$ and $\xi$ points of spectral concentration.
Thus we have largely confined our search for the transitional behaviour of $\theta(x,\mu)$ to the $x-$range $(0,18)$.
 There remains the question, on which we have touched in the remarks relating
to Figure 3, whether $\theta(x,\mu)$ possesses the transitional property also for higher values of $r$ and, more particularly, whether an infinity of values of $r$ is involved in this way for some fixed value of $c$.
Thus we have the theoretical question of whether an infinite set (bounded or unbounded) of spectral concentration points can exist for some $c$.
All that is known is that an unbounded set cannot occur when $ a >2$ because it is shown in \cite{MSPE98} that the set of spectral concentration points is bounded for any $q$ in (\ref{eq:1.2}) such that $x q(x)$ is $L(0,\infty)$.

   \bibliographystyle{plain}
\bibliography{../bib_dir/bib1,../bib_dir/bibliography,../bib_dir/help,../bib_dir/specon,../bib_dir/specon2}

\end{document}